\documentclass{cimart}

\usepackage{tikz} \usetikzlibrary{patterns, positioning, calc}
\usepackage[normalem]{ulem}

\newcommand{\E}{{\mathbb E}}
\newcommand{\eps}{{\varepsilon}}
\newcommand{\R}{{\mathbb R}}

\def\vint_#1{\mathchoice%
	{\mathop{\kern 0.2em\vrule width 0.6em height 0.69678ex depth -0.58065ex
			\kern -0.8em \intop}\nolimits_{\kern -0.4em#1}}%
	{\mathop{\kern 0.1em\vrule width 0.5em height 0.69678ex depth -0.60387ex
			\kern -0.6em \intop}\nolimits_{#1}}%
	{\mathop{\kern 0.1em\vrule width 0.5em height 0.69678ex
			depth -0.60387ex
			\kern -0.6em \intop}\nolimits_{#1}}%
	{\mathop{\kern 0.1em\vrule width 0.5em height 0.69678ex depth -0.60387ex
			\kern -0.6em \intop}\nolimits_{#1}}}

\definecolor{darkblue}{rgb}{0.05, .05, .9}
\definecolor{darkgreen}{rgb}{0.1, .65, .1}
\definecolor{darkred}{rgb}{0.8,0,0}

\title[A two-player zero-sum probabilistic game that approximates the mean curvature flow]{A two-player zero-sum probabilistic game that approxi\\mates the mean curvature flow}

\author{Irene Gonz\' alvez, Alfredo Miranda, Julio D. Rossi and Jorge Ruiz-Cases}

\authorinfo[I. Gonz\'alvez]{Department of Mathematics, Autonomous University of Madrid, Spain}{irene.gonzalvez@uam.es}

\authorinfo[A. Miranda]{Department of Mathematics, University of Buenos Aires, Argentina}{amiranda@dm.uba.ar}

\authorinfo[J. D. Rossi]{Department of Mathematics and Statistics, Torcuato Di Tella University, Argentina}{julio.rossi@utdt.edu}

\authorinfo[J. Ruiz-Cases]{Department of Mathematics, Autonomous University of Madrid, Spain}{jorge.ruizc@uam.es}

\abstract{In this paper we introduce a new two-player zero-sum game whose value function approximates the level set formulation for the geometric evolution by mean curvature of a hypersurface. In our approach the game is played with symmetric rules for the two players and probability theory is involved (the game is not deterministic).}



\VOLUME{34}
\YEAR{2026}
\ISSUE{2}
\NUMBER{1}
\DOI{https://doi.org/10.46298/cm.15764}

\begin{document}

\section{Introduction}

\subsection{Description of the main goal}

The aim of this paper is to study a two-players zero-sum probabilistic game 
whose value functions approximates the motion by mean curvature of a hypersurface. 
We show that value functions of the game converge, as a parameter that controls
the size of the steps made in the game goes to zero, to the unique viscosity solution
to a partial differential equation whose level sets encode the motion by mean curvature.

Game theory was used to obtain a better understanding of the mean curvature flow in 
\cite{KS} and \cite{Misu}. In those references, the authors considered a deterministic (no probability theory is involved) two-person zero-sum game 
with asymmetric rules for the players. 
Our goal here is to introduce a different game in which the two players play
with symmetric rules and some randomness is involved.
In order to pass to the limit in the value functions of the game and find mean curvature as the limit equation we use viscosity solutions (we will rely on previous theory developed for viscosity solutions to geometric flows, see \cite{CIL} and \cite{GGIS}).

Motion of a hypersurface by mean curvature is by now well-understood. Its usual interpretation involves the steepest-descent for
the perimeter functional. Here we will use a different perspective (using game theory) to 
study this geometric movement.
Convexity is preserved under motion by mean curvature, so the boundary of
a convex body shrinks monotonically.
Thus, when the initial hypersurface is the boundary of a convex domain, the
mean curvature flow can be described in two equivalent ways: by
following the moving boundary as an evolving surface, or by specifying for each point inside the domain the arrival
time when the moving hypersurface reaches the point. 
We follow the second interpretation and in our game-theoretic interpretation, the latter viewpoint is
associated with a minimum-exit-time problem.

Now, let us describe the main ingredients that we need to state and prove our results. First we introduce an
elliptic nonlinear partial differential equation (PDE) related to the geometric motion by mean curvature and 
next we describe a game whose value functions
approximate solutions to the PDE.

\subsubsection{The movement by mean curvature of a hypersurface
and its associated elliptic equation}

Our aim is to describe how a hypersurface that is the boundary of a connected and
strictly convex domain, $S=\partial \Omega_{0}\subset \mathbb{R}^N$,  $N\geq 2$, evolves  
according to the mean curvature flow. We will use a level set
approach to describe this geometric evolution. Assume that there is a 
real valued function, $u(x)$, defined for $x \in \Omega_0$, and consider
the $t$ superlevel sets of $u(x)$, 
$$
\Omega_t = \{ x : u (x) > t \},\quad t\geq 0.
$$
In what follows we denote by $\nabla u (x)$
and by $D^2 u (x)$ the gradient and the Hessian of $u$ 
with respect to the spacial variable, $x$. 
Assume that $\partial \Omega_t$ is smooth. 
Take $x\in \partial \Omega_t$  a regular point ($x$ is a regular point if  
$\nabla u (x)\neq 0$). We have 
$\nabla u (x) \perp \partial \Omega_t$ and 
for a unitary vector $v \perp \nabla u (x)$ (notice that $v$ is
tangential to the hypersurface $\partial \Omega_t$)
the quantity $-\langle D^2 u (x) v , v \rangle$ gives the  
curvature of $\partial \Omega_t$ in the direction of $v$.  
Therefore, under these conditions,
the mean of the principal curvatures of $\partial \Omega_t$ at a regular point is given by
$$
\begin{array}{l}
\displaystyle \kappa = \sum_i \kappa_i =  \mbox{div} \Big(\frac{\nabla u}{|\nabla u|}\Big)(x) =\frac{1 }{|\nabla u (x) |}
\Big(\Delta u (x)  - \Big\langle D^2 u (x) \frac{\nabla u}{|\nabla u|} (x),  \frac{\nabla u}{|\nabla u|} (x) \Big\rangle \Big).
\end{array}
$$
We consider the geometric evolution of the hypersurface
$\partial \Omega_t$ moving its points in the direction of the normal vector (pointing
inside the set $\Omega_t$) with speed given by the mean curvature,  
$V = -\kappa$ on $\partial \Omega_t$. 
Now, we look at the elliptic equation 
\begin{equation}\label{1} 
    \mathcal{L}(u(x)) := \Delta u (x) - \Big\langle D^2 u(x) \frac{\nabla u}{|\nabla u|}(x) ,  \frac{\nabla u}{|\nabla u|}(x) \Big\rangle =-1,
\end{equation}
and its associated Dirichlet problem,
	\begin{equation}
		\label{PDEPb}
\begin{cases}
	 \mathcal{L}(u(x)) = - 1,& \;\; x\in
	  \Omega_0, \\
	u(x) =0,&\;\;x\in\partial \Omega_0.
	\end{cases}
   \end{equation}

For this problem it is known that there exists a unique viscosity solution
and in addition a comparison principle holds, see \cite{ES,KS} and references therein.
The analysis of \eqref{PDEPb} in \cite{ES,KS} uses the framework of viscosity solutions. This is necessary
because in its classical form the PDE  is not well-defined when
$\nabla u =0$. However, the problem has nice solutions. Indeed, for a convex planar domain, the  
evolution remains smooth under motion
by mean curvature, and it becomes asymptotically circular as it shrinks to a point \cite{10}. Concerning
regularity, in \cite{KS} it is proved that
$u$ is $C^3 (\Omega_0)$.

Notice that when $u$ solves \eqref{PDEPb} then $v(x,t) = u(x) -t$ is a solution to 
$v_t (x,t) =  \mathcal{L}(v(x,t))$ that is usually referred to as the mean curvature evolution equation. 
As we have mentioned, for each point inside the domain, we look for the arrival
time of the hypersurface. Remark that $
\Omega_t = \{ x : v (x,t) > 0 \} = \{ x : u (x) > t \}$. 
Therefore, \eqref{PDEPb} is the PDE associated to the level set formulation
for the first arrival time for the geometric evolution of a convex hypersurface by mean curvature.

The motion by mean curvature is nowadays a classical subject, we refer to~\cite{2,ES,CGG,Giga,8,15,16} and also to~\cite{9,Mercier,19,Misu}. 
Concerning game theoretical approximations for this problem, we quote
\cite{KS}, where the authors introduce a two-players zero-sum game whose value functions approximate solutions to \eqref{PDEPb}. This game 
is played with asymmetric rules among the players (one chooses a direction and 
the other a sign) and does not
involve probability (the next position of the game only depends 
on the choices made by the players). See also \cite{Misu,Nosotros}  for other variants of this game. 
As we have mentioned, here we introduce a new game in which the two players play
with symmetric rules and some randomness is involved in the
choice of the next position of the game.

\subsubsection{A probabilistic game approximation for the elliptic problem}
Next, let us describe a game whose value function approximates 
the solution to \eqref{PDEPb}.

The game is a probabilistic two-person zero-sum game. 
As in \cite{KS} we use Paul for the name of the first player and Carol for the second player. Take $\varepsilon>0$ (a parameter that controls the size of the possible movements in the game), $\Omega_0 \subset \mathbb{R}^N$ a 
strictly convex
and bounded domain. 
Let $x_0\in \Omega_0$ be the initial position of the game. The game is played as follows: at the $i-$th round, Paul chooses a set of unitary vectors
$A_i\subset S^{N-1}$ with surface measure $\sigma(A_i) \geq 
\frac12 \sigma (S^{N-1}) + \delta_\varepsilon$ (with $\delta_\varepsilon \approx \eps^{1/2}$) and Carol chooses a set
$B_i\subset S^{N-1}$ with measure $\sigma(B_i) \geq 
\frac12 \sigma (S^{N-1}) + \delta_\varepsilon$. Here and in what follows we 
denote by $\sigma$ the surface measure on the sphere $S^{N-1}$.
Notice that both players choose a set slightly bigger than half of the sphere
and therefore the intersection of both sets, $A_i\cap B_i$, has positive measure.
Once these choices are made, the next position of the game 
is given by
\[x_i = x_{i-1} + v_i \varepsilon,\]
where the vector $v_i$ is randomly chosen (with uniform probability)
in the set $A_i\cap B_i$.
The game ends when the position exits $\Omega_0$ and Carol pays to Paul
an amount proportional to the number of plays, that is, the payoff is given by 
$\varepsilon^2 K \times (\mbox{number of plays})$ with $K$ a constant
that we will specify latter.

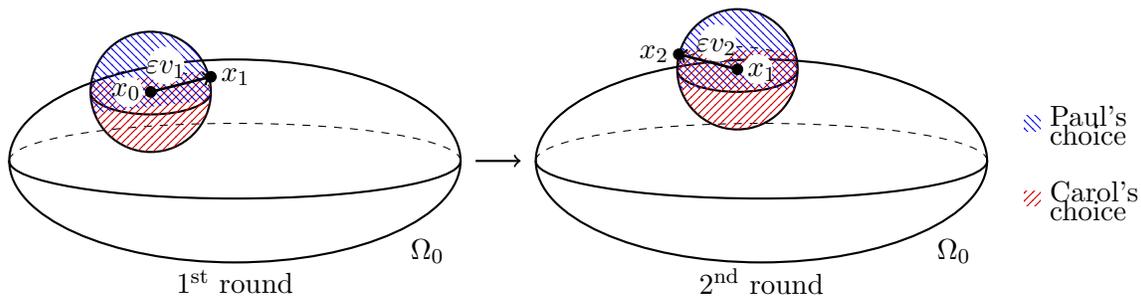
\begin{figure}[h!]
	\centering
	\begin{tikzpicture}[scale=1.0]
		%\tableros
		\draw[thick] (0,0) ellipse (3cm and 1.35cm); %centro
		\begin{scope}%tercero
			\clip (-3.1,-1.6) rectangle (3.1,0); % Pinta la mitad inferior 
			\draw[thick] (0,0) ellipse (3cm and 0.5cm);
		\end{scope}
		\begin{scope} % Pinta la mitad superior
			\clip (-3.1,1.6) rectangle (3.1,0);  
			\draw[dashed] (0,0) ellipse (3cm and 0.5cm);
		\end{scope}
		\draw[thick] (-7,0) ellipse (3cm and 1.35cm); %izq
		\begin{scope}%tercero
			\clip (-3.1-7,-1.6) rectangle (3.1-7,0); % Pinta la mitad inferior 
			\draw[thick] (-7,0) ellipse (3cm and 0.5cm);
		\end{scope}
		\begin{scope} % Pinta la mitad superior
			\clip (-3.1-7,1.6) rectangle (3.1-7,0);  
			\draw[dashed] (-7,0) ellipse (3cm and 0.5cm);
		\end{scope}

		%transicion
		\draw[thick,->](-3.8,0)--(-3.2,0);
		
		%PrimeraRonda-Primer Tablero 
		
		\draw[thick] (-1.3*0.707-0.2-7,1.3*0.707) circle (0.8cm); %S^N
		\begin{scope}
			\clip (-1.3*0.707-0.2-7-0.8,1.3*0.707-0.4) rectangle (-1.3*0.707-0.2-7+0.8,1.3*0.707); % Pinta la mitad inferior 
			\draw[thick] (-1.3*0.707-0.2-7,1.3*0.707)  ellipse (0.8cm and 0.3cm);
		\end{scope}
		\begin{scope}
			\clip (-1.3*0.707-0.2-7-0.8,1.3*0.707+0.4) rectangle (-1.3*0.707-0.2-7+0.8,1.3*0.707); % Pinta la mitad superior
			\draw[dashed] (-1.3*0.707-0.2-7,1.3*0.707)  ellipse (0.8cm and 0.3cm);
		\end{scope}
		%Primera Ronda- Elecciones
		%Eleccion de Carol
		% \draw[dashed](-1.3*0.707-0.2-7-0.9,1.3*0.707-0.3) --(-1.3*0.707-0.2-7+0.9,1.3*0.707-0.1);
		%	\draw[dashed](-1.3*0.707-0.2-7-0.9,1.3*0.707+0.1) --(-1.3*0.707-0.2-7+0.9,1.3*0.707+0.3);
		\begin{scope}
			\clip(-1.3*0.707-0.2-7-0.9,1.3*0.707+0.1) --(-1.3*0.707-0.2-7+0.9,1.3*0.707+0.3) --(-1.3*0.707-0.2-7+0.9,1.3*0.707-0.1-0.8)--(-1.3*0.707-0.2-7-0.9,1.3*0.707-0.3-0.8) -- cycle; % Pinta la mitad inferior
			\fill[pattern=north east lines, pattern color=darkred]  (-1.3*0.707-0.2-7,1.3*0.707) circle (0.8cm);
		\end{scope}
		%Eleccion de Paul
		\begin{scope}
			\clip(-1.3*0.707-0.2-7-0.9,1.3*0.707-0.3) --(-1.3*0.707-0.2-7+0.9,1.3*0.707-0.1) --(-1.3*0.707-0.2-7+0.9,1.3*0.707-0.1+1)--(-1.3*0.707-0.2-7-0.9,1.3*0.707-0.3+1.8) -- cycle; % Pinta la superior
			\fill[pattern=north west lines, pattern color=darkblue] (-1.3*0.707-0.2-7,1.3*0.707) circle (0.8cm);
		\end{scope}
		%Vector v
		\draw[thick, ->, line width=1pt](-1.3*0.707-0.2-7,1.3*0.707) --(-1.3*0.707-0.2-7+0.8,1.3*0.707+0.2);
		%puntos
		\filldraw[black] (-1.3*0.707-0.2-7,1.3*0.707) circle (2pt); %x_0
		\fill[white ] (-1.3*0.707-0.2-7.3,1.3*0.707) circle (0.2cm);
		\node[left] at (-1.3*0.707-0.2-7,1.3*0.707) {$x_0$};
		\fill[white ] (-1.3*0.707-0.2-6.8,1.3*0.707+0.35) circle (0.25cm);
		\node[left] at (-1.3*0.707-0.2-6.4,1.3*0.707+0.3) {$\varepsilon v_1$};
		\filldraw[black] (-1.3*0.707-0.2-7+0.8,1.3*0.707+0.2) circle (2pt); %x_1
		\node[right] at (-1.3*0.707-0.2-7+0.8,1.3*0.707+0.2) {$x_1$};
		\node[right] at (2.2-7,-1.2) {$\Omega_0$};
		\node at (-7,-1.6) {$1^{\textrm{st}}$ round};
		%%%%%%%%%%%%%%%%%%%%
		%%%%%%%%%%%%%%%%%%%%%
		%%%%%%%%%%%%%%%%%%%%%%%   
		%SegundaRonda-SegundaTablero 
		\draw[thick] (-1.3*0.707-0.2+0.8,1.3*0.707+0.3)  circle (0.8cm); %S^N
		\begin{scope}
			\clip (-1.3*0.707-0.2+0.8-0.8,1.3*0.707+0.3-0.4) rectangle (-1.3*0.707-0.2+0.8+0.8,1.3*0.707+0.3); % Pinta la mitad inferior 
			\draw[thick] (-1.3*0.707-0.2+0.8,1.3*0.707+0.3)  ellipse (0.8cm and 0.3cm);
		\end{scope}
		\begin{scope}
			\clip (-1.3*0.707-0.2+0.8-0.8,1.3*0.707+0.3) rectangle (-1.3*0.707-0.2+0.8+0.8,1.3*0.707+0.3+0.3); % Pinta la mitad superior
			\draw[dashed] (-1.3*0.707-0.2+0.8,1.3*0.707+0.3)  ellipse (0.8cm and 0.3cm);
		\end{scope}
		%Eleccion de Paul
		\begin{scope}
			\clip(-1.3*0.707-0.2+0.8-0.8,1.3*0.707+0.3-0.25) --(-1.3*0.707-0.2+0.8+0.8,1.3*0.707+0.3-0.25) --(-1.3*0.707-0.2+0.8+0.8,1.3*0.707+0.3+0.9)--(-1.3*0.707-0.2+0.8-0.8,1.3*0.707+0.3+0.9) -- cycle; % Pinta la mitad inferior
			\fill[pattern=north west lines, pattern color=darkblue]  (-1.3*0.707-0.2+0.8,1.3*0.707+0.3) circle (0.8cm);
		\end{scope}
		%Eleccion de Carol 
		\begin{scope}
			\clip(-1.3*0.707-0.2+0.8-0.8,1.3*0.707+0.3+0.25) --(-1.3*0.707-0.2+0.8+0.8,1.3*0.707+0.3+0.25) --(-1.3*0.707-0.2+0.8+0.8,1.3*0.707+0.3-0.9)--(-1.3*0.707-0.2+0.8-0.8,1.3*0.707+0.3-0.9) -- cycle; % Pinta la mitad inferior
			\fill[pattern=north east lines, pattern color=darkred]  (-1.3*0.707-0.2+0.8,1.3*0.707+0.3) circle (0.8cm);
		\end{scope}
		%Vector v
		\draw[thick, ->, line width=1pt](-1.3*0.707-0.2+0.8,1.3*0.707+0.3) --(-1.3*0.707-0.2+0.8-0.78,1.3*0.707+0.3+0.2);
		%Puntos 
		\fill[white ] (-1.3*0.707-0.2+0.8+0.27,1.3*0.707+0.3) circle (0.2cm);
		\node[right] at (-1.3*0.707-0.2+0.8,1.3*0.707+0.26) {$x_1$};
		\filldraw[black] (-1.3*0.707-0.2+0.8,1.3*0.707+0.3) circle (2pt); %x_1
		\fill[white ] (-1.3*0.707-0.2+0.8-0.78+0.52,1.3*0.707+0.3+0.2+0.15
		) circle (0.25cm);
		\node[left] at (-1.3*0.707-0.2+0.8-0.78+0.9,1.3*0.707+0.3+0.2+0.1) {$\varepsilon v_2$};
		\filldraw[black] (-1.3*0.707-0.2+0.8-0.78,1.3*0.707+0.3+0.2) circle (2pt); %x_2
		\node[left] at (-1.3*0.707-0.2+0.8-0.78,1.3*0.707+0.3+0.2) {$x_2$};
		\node[right] at (2.2,-1.2) {$\Omega_0$};
		\node at (0,-1.6) {$2^\textrm{nd}$ round};
		%Leyenda
		\node[right] at (3.7,0.6) {Paul's};
		\node[right] at (3.7,0.35) {choice};
		\fill[pattern=north west lines, pattern color=darkblue](3.7-0.1-0.1,0.5+0.1) -- (3.7+0.1-0.1,0.5+0.1) -- (3.7+0.1-0.1,0.5-0.1) -- (3.7-0.1-0.1,0.5-0.1)-- cycle;
		\node[right] at (3.7,-0.4) {Carol's};
			\node[right] at (3.7,-0.65) {choice};
		\fill[pattern=north east lines, pattern color=darkred] (3.7-0.1-0.1,-0.5+0.1)--(3.7+0.1-0.1,-0.5+0.1)--(3.7+0.1-0.1,-0.5-0.1)--(3.7-0.1-0.1,-0.5-0.1)-- cycle; 
	\end{tikzpicture}
		\caption{ \label{FigureGame} 
			Game starting at $x_0$.}
\end{figure}

We define $u^\varepsilon (x)$ as the value of this game starting at $x$ (the value of the game is just the expected final payoff optimized by both players, Paul wants to maximize the expected outcome, while Carol aims to minimize it). The value of the game for Paul is given by the following formula, 
\[u_p^\varepsilon (x_0) \! = \inf_{S_c} \sup_{S_p} \mathbb{E}^{x_0}_{S_p,S_c} 
\Big[\varepsilon^2 K \times (\mbox{number of plays})\Big].
\]
Here the inf and sup are taken among the possible strategies for Paul and Carol, that we denote by $S_p$ and $S_c$,
respectively.
This is the best possible outcome that Paul may obtain, provided that both
players play their best.
Analogously, the value for Carol is given by
\[u_c^\varepsilon (x_0) \! =  \sup_{S_p} \inf_{S_c} \mathbb{E}^{x_0}_{S_p,S_c} 
\Big[ \varepsilon^2 K \times (\mbox{number of plays}) \Big].
\]
Notice that we always have $u_p (x_0) \geq u_c(x_0)$ for every $x_0 \in \Omega$.
Finally, the value of the game is defined as 
$$
u^\varepsilon (x_0) := u_p^\varepsilon (x_0) = u_c^\varepsilon (x_0) 
$$
provided the values for Paul and Carol coincide.

For our game there is a value that is given by the unique solution to 
\begin{equation}
	\label{DPP} 
		\left\{
	\begin{array}{ll}
	\displaystyle
	u^\varepsilon(x) = \sup_{A} \inf_{B}  \left\{ 
\vint_{A\cap B}  u^\varepsilon \big(x + v\varepsilon \big) 
d \sigma (v)\right\} +  \varepsilon^2 K,  & x\in \Omega_0,\\[10pt]
	u^\varepsilon(x) = 0 ,  & x\in \mathbb{R}^N \setminus 
\Omega_0.
	\end{array} \right.
\end{equation}
To simplify the notation, each time we compute a supremum or an infimum among sets we 
understand that we are considering subsets of $S^{N-1}$ with surface measure
bigger or equal than $ 
\frac12 \sigma (S^{N-1}) + \delta_\eps$. 
We use $\vint_C f(v) d \sigma (v)$ for the average
$$
\vint_C f(v) d \sigma (v) :=  \frac{1}{\sigma(C)}
\int_{C}  f(v)
d \sigma (v).
$$ 

The equation that appears in \eqref{DPP} is known as the Dynamic Programming Principle (DPP) in the game literature,
see \cite{MS2}, and reflects the rules of the game. In fact, 
the equation that appears in \eqref{DPP} shows the outcome after only one round of the game starting from $x$.
Paul wants to
maximize the outcome, and he can choose the set $A$,
while Carol aims to minimize, and chooses the set $B$. Then the outcome
or the value after making one move is given by the supremum among Paul's choices of the minimum 
among Carol's choices of the mean value at the new position 
(the new position is random inside $A\cap B$) plus $K \eps^2$ (the game counts 
one extra move). 
We observe that Paul plays aiming to stay inside $\Omega_0$ for a large number of plays
while Carol tries to exit $\Omega_0$ as quickly as possible.

Now, to look for a partial differential equation related to this game, we argue formally (one of the main goals of this paper is to make rigorous what follows). 
Let us find the asymptotic behaviour as $\varepsilon \approx 0$ of the 
DPP, \eqref{DPP}, evaluated at a smooth function. Assume that for a $C^{2,1}$ function $\phi$ we have
\begin{equation}
	\label{DPP-2233} \phi (x) \approx 
 \sup_{A} \inf_{B}  \left\{ \frac{1}{\sigma(A\cap B)}
\int_{A\cap B}  \phi \big(x + v\varepsilon \big) 
d \sigma (v)\right\} +  \varepsilon^2K, \end{equation}
that is, 
\begin{equation}
	\label{DPP-33} 0 \! \approx \!   
 \sup_{A} \inf_{B}  \left\{ \frac{1}{\sigma(A\cap B)}
\int_{A\cap B} ( \phi \big(x + v\varepsilon \big)  - \phi (x)  )
d \sigma (v)\right\} +  \varepsilon^2 K.
\end{equation}
Using that $\phi \in C^{2}$ and neglecting higher order terms, from a simple Taylor expansion, we arrive to 
\begin{equation}
	\label{DPP-44} 0 \! \approx \!   
 \sup_{A} \inf_{B}  \left\{ \frac{1}{\sigma(A\cap B)}
\int_{A\cap B} \frac{1}{\varepsilon} \langle \nabla \phi (x), v \rangle +
\frac12
\langle D^2 \phi (x) v,v \rangle )
d \sigma (v)\right\} + K.
\end{equation}

Now, assuming that $\nabla \phi (x) \neq 0$, the leading term is the one that involves $1/\varepsilon$.
Hence, when computing the supremum Paul wants to use a vector $v$ such that
 $ \langle \nabla \phi (x), v \rangle $ is as large as possible, and 
 he may choose $A= \{ v:  \langle \frac{\nabla\phi(x)}{|\nabla \phi (x)|} , v \rangle \geq 
-\theta_\varepsilon \}$ with $\theta_\eps$ such that $\sigma(A) = 
\frac12 \sigma (S^{N-1}) + \delta_\varepsilon$ (remark that we have the constraint $\sigma(A) \geq 
\frac12 \sigma (S^{N-1}) + \delta_\varepsilon$). Notice that since $\delta_\eps \approx \eps^{1/2}$
we also have $\theta_\eps \approx \eps^{1/2}$, in fact, $\sigma (A) \approx 
\frac12 \sigma (S^{N-1}) + C \theta_\varepsilon $ for $\theta_\varepsilon$ small. 
Analogously, a  clever choice for Carol could be $B= \{ v:  \langle \frac{\nabla\phi(x)}{|\nabla \phi (x)|}, v \rangle \leq 
\theta_\varepsilon \}$. With these choices of the sets $A$ and $B$ we get
$A\cap B = \{  v : -\theta_\varepsilon \leq
\langle \frac{\nabla\phi(x)}{|\nabla \phi (x)|}, v \rangle \leq 
\theta_\varepsilon\}$ and, in this case, by symmetry, the leading  term 
vanishes, since we have
$$
\int_{-\theta_\varepsilon \leq
\langle \frac{\nabla\phi(x)}{|\nabla \phi (x)|}, v \rangle \leq 
\theta_\varepsilon}  \langle \nabla \phi (x), v \rangle d \sigma (v)
=0.
$$
Then, from \eqref{DPP-44} with the sets $A$ and $B$ described before, we arrive to 
\begin{equation}
	\label{DPP-44-66} 0 \! \approx \!   
\frac{1}{\sigma ( \{ v: -\theta_\varepsilon \leq
\langle \frac{\nabla\phi(x)}{|\nabla \phi (x)|}, v \rangle \leq 
\theta_\varepsilon \} )}\int_{-\theta_\varepsilon \leq
\langle \frac{\nabla\phi(x)}{|\nabla \phi (x)|}, v \rangle \leq 
\theta_\varepsilon}\frac12
\langle D^2 \phi (x) v,v \rangle )
d \sigma (v) + K.
\end{equation}
Observe that the sets $A\cap B = \{v: \, -\theta_\varepsilon \leq
\langle \frac{\nabla\phi(x)}{|\nabla \phi (x)|}, v \rangle \leq 
\theta_\varepsilon \}$ converge to $\{v:\langle \nabla \phi (x), v \rangle=0\}$
as $\eps\to 0$ 
and hence, passing to the limit as $\eps \to 0$, we obtain 
\begin{equation}
	\label{DPP-55} 0 \! = \!   
 \frac{1}{\mu (\langle \nabla \phi (x), v \rangle=0)}
\int_{\langle \nabla \phi (x), v \rangle=0} \frac12
\langle D^2 \phi (x) v,v \rangle )
d \mu (v) + K.
\end{equation}
Notice that here $\mu$ stands for the surface measure in dimension $N-2$
(we integrate on the set $\{ v \in S^{N-1}: \langle \nabla \phi (x), v \rangle=0\} = S^{N-2}$).
At this point we highlight that a rigorous proof of this whole limit procedure is one
of the main difficulties in this article. 

Now, we just observe that we can compute the Laplacian of $ \phi$ at $x$
using an orthonormal base in which one of the vectors is in the direction
of the gradient. We have, assuming for simplicity that $ \nabla \phi (x)$  points
in the direction of $v_N$,
$$
\begin{array}{l}
\displaystyle
 \frac{1}{\mu(\langle \nabla \phi (x), v \rangle=0)}
\int_{\langle \nabla \phi (x), v \rangle=0} \frac12
\langle D^2 \phi (x) v,v \rangle )
d \mu (v) \\[10pt]
\displaystyle
\qquad =\sum_{i,j=1}^{N-1} \frac{1}{\mu(\langle \nabla \phi (x), v \rangle=0)}
\int_{\langle \nabla \phi (x), v \rangle=0} \frac12
\langle D^2 \phi (x) v_i,v_j \rangle ) d \mu (v) \\[10pt]
\displaystyle
\qquad  = C \sum_{i=1}^{N-1}
\frac{\partial^2 \phi}{\partial v_i^2} (x) = C
\Delta \phi (x) \mid_{\langle \nabla \phi (x), v \rangle=0}
\end{array}
$$
for a constant $C$ given by 
$$
C= \frac{1}{2\mu(\{\langle e_N, v \rangle=0\})} \int_{\langle e_N, v \rangle=0} (v_1)^2 \, d \mu (v).
$$
Then, \eqref{DPP-55} is just
\begin{equation}
	\label{DPP-66} 0 \! = \!   
C
\Delta \phi (x) \mid_{\langle \nabla \phi (x), v \rangle=0}
 + K.
\end{equation}
Hence, if we choose $K=C$, we arrive to the mean curvature equation 
\begin{equation}
	\label{DPP-77} 0 \! = \!   
\Delta \phi (x) \mid_{\langle \nabla \phi (x), v \rangle=0}
 + 1 =  \Delta \phi (x) - \Big\langle D^2 \phi (x) \frac{\nabla \phi}{|\nabla \phi|},  \frac{\nabla \phi }{|\nabla \phi|}
 \Big\rangle + 1.
\end{equation}
For other choices of $K>0$ we find the equation associated with the movement by a scalar multiple
of the mean curvature.

Therefore, we conclude that the DPP associated to the game, \eqref{DPP}, is related to the level set formulation of the mean curvature flow. 
Then, we expect that the limit as $\varepsilon \to 0$ for the value function of the game converges to 
the viscosity solution to \eqref{PDEPb}.

To make this argument rigorous we run into several difficulties, 
First, to use viscosity arguments to test the equation with smooth test
functions we need to handle inequalities. When we face an inequality 
we can fix one of the sets, $A$ or $B$, as above, and still continue 
with an inequality, however, in order to pass to the limit, we need estimates that hold
regardless the choice of the other set. This fact forces us to prove a delicate
technical lemma from geometric measure theory, see Lemma \ref{lema.clave}. 
 On top of this, in the limit we need to obtain the boundary condition $u=0$ on 
$\partial \Omega_0$. Since we always have that the values of the game 
$u^\varepsilon$ are positive inside $\Omega_0$ (the game is played at least
one time before exiting the domain) we only need an upper bound for 
$u^\varepsilon (x)$ for points $x$ inside $\Omega_0$ that are close to the boundary. 
The fact that this bound has some uniformity in $\varepsilon$ 
(this is needed since we want to pass to the limit as $\varepsilon \to0$)
relies on the fact that we assumed that 
 $\Omega_0$ is strictly convex (this assumption also appears in
\cite{KS}).

Finally, we introduce the sets
$$
\Omega_t^{\varepsilon} = \{ x: u^\varepsilon (x) >t \},
$$
that is, the set of positions in $\mathbb{R}^N$ where the players expect to play more than  $\left\lceil \frac{t}{k \varepsilon^2} \right\rceil$ rounds.
As a consequence of the previous locally uniform convergence
of $u^\eps$ to $u$, we have that for each $t>0$, when $\varepsilon$ is small the set where the value function of the $\varepsilon$-game is greater than $t$, $\Omega_t^\eps = \{x : u^\eps (x)>t\}$ is close to the set where the solution to the mean curvature problem is greater than $t$, $\Omega_t = \{x : u (x)>t\}$.

\subsection{Statements of the main results}

Let us state rigorously the main results that are included in this
paper. 
First, let us describe our results for the game.
 In Section \ref{sect.Game} we include the following results.

\begin{theorem}\label{maintheorem.DPP}
There exists a unique solution to \eqref{DPP}. 
Moreover, a comparison principle between super and subsolutions to \eqref{DPP} holds. 
\end{theorem}

With this result at hand we can show that the game has a value that coincides with 
the solution to \eqref{DPP}.

\begin{theorem}\label{maintheorem.DPP-33}
The game has a value $u^{\varepsilon}$  that is characterized as the unique solution to \eqref{DPP}. 
\end{theorem}

Next,  in Section \ref{sect.GameAndEq} we deal with the limit as $\varepsilon \to 0$
of the value function of the game, $u^{\varepsilon}$. Our main result in this paper is
the following: 

\begin{theorem}\label{maintheorem.juego.conv}
 Let $\{u^{\varepsilon}\}_{\varepsilon>0}$  be a family of value functions of the game 
 inside $\Omega_0$,
	that is, for each $\varepsilon>0$,  $u^{\varepsilon}$ is the solution to \eqref{DPP}. 
	
	Then, $u^{\varepsilon}$ converges uniformly in $\overline{\Omega}_0$ to $u$, the unique solution to \eqref{PDEPb},
	as $\varepsilon \to 0$.
\end{theorem}

As a consequence, we can obtain the behaviour of
the positivity sets of $u^\varepsilon$ as $\varepsilon \to 0$.

\begin{corollary}\label{cor1} 
Consider  
	\begin{equation}
	\label{OmegaSets}	\Omega^{\varepsilon}_{t}:=\{x\,:\;\;u^{\varepsilon}(x)>t\}\;\;\textrm{and}\;\;\Omega_{t}:=\{x\,:\;\;u(x)>t\}.
	\end{equation}Then, for each $t>0$, we have that 
	\begin{equation*}
		\Omega_t\subset \liminf_{\varepsilon \to 0} \Omega^\varepsilon_t \subset \limsup_{\varepsilon \to 0} \Omega^\varepsilon_t\subset \overline{\Omega}_t.
	\end{equation*}
\end{corollary}

\subsection{Related results} 
Let us end this introduction with extra references of related results in the literature. 

The level set method for the study of evolution equations for
hypersurfaces was first rigorously analyzed in 
\cite{CGG,ES}. This method was applied to various equations including the mean curvature flow equation
(with or without an obstacle), see \cite{Almeida,CGG,ES}. For general references on 
level sets formulation of geometric flows we refer to \cite{2,6,8,9,15,16,19} and the book \cite{Giga}.

On the other hand, the relation between game theory and nonlinear PDEs 
is quite rich and has attracted considerable attention in recent years. We quote \cite{BR,KS,KS2,Luiro,MPR,MPR2,ML,MiRo,Misu,PSSW,PS,R} and the books \cite{BRLibro,Lewicka}.
Concerning games for geometric flows, as we have mentioned, our immediate precedent is \cite{KS} where
the authors analyze a deterministic game related to the mean curvature flow.
For the formulation of this game in the presence of an obstacle we refer to  \cite{Misu}.

\section{Preliminaries}\label{sect.Prelim}

\subsection{The mean curvature equation}

First, we collect some information for the mean curvature problem \eqref{PDEPb}.

\subsubsection{Viscosity solutions}
Following \cite{KS}, let us state the precise definition of being a viscosity solution
to our elliptic problem,
	\begin{equation}
		\label{PDEPb-2-E}
\begin{cases}
\displaystyle	\Delta u(x) - \Big\langle   D^2  u (x) \frac{\nabla u}{|\nabla u|} (x), 
 \frac{\nabla u}{|\nabla u|} (x)
\Big \rangle  = - 1,& \;\; x\in
	  \Omega_0, \\
	u(x) =0,&\;\;x\in\partial \Omega_0.
	\end{cases}
   \end{equation}

\begin{definition} \label{def-1}
An uppersemicontinuous function $\underline{u}$ is a subsolution to \eqref{PDEPb-2-E}
provided that if
$\underline{u} - \phi$
 has a local maximum at a point $x_0 \in \Omega_0$
for some $\phi \in C^2$, then
$$
\Delta \phi (x_0) - \Big \langle   D^2  \phi (x_0) \frac{\nabla \phi}{|\nabla \phi|} (x_0), 
 \frac{\nabla \phi}{|\nabla \phi|} (x_0)
\Big\rangle \geq - 1
$$
when $\nabla \phi (x_0) \neq 0$; and 
$$
\Delta \phi (x_0) - \Big\langle  D^2  \phi (x_0) \eta, \eta 
 \Big\rangle  \geq -1,
$$
for some vector $\eta$ with $|\eta|\leq 1$;
and $\underline{u} \leq 0$ on $\partial \Omega$.

A lowersemicontinuous function $\overline{u}$ is a supersolution to \eqref{PDEPb-2-E}
provided that if
$\overline{u} - \phi$
 has a local minimum at a point $x_0 \in \Omega_0$
for some $\phi \in C^2$, then
$$
\Delta \phi (x_0) - \Big\langle   D^2  \phi (x_0) \frac{\nabla \phi}{|\nabla \phi|} (x_0), 
 \frac{\nabla \phi}{|\nabla \phi|} (x_0)
\Big\rangle \leq - 1
$$
when $\nabla \phi (x_0) \neq 0$; and 
$$
\Delta \phi (x_0) - \Big\langle  D^2  \phi (x_0) \eta, \eta 
\Big\rangle   \leq -1,
$$
for some vector $\eta$ with $|\eta|\leq 1$;
and $\overline{u} \geq 0$ on $\partial \Omega$.

Finally, a function $u$ is a solution to \eqref{PDEPb-2} if it is both a
super and a subsolution.
\end{definition}

We can rewrite \eqref{PDEPb-2-E} as 	
\begin{equation}
		\label{PDEPb-2}
\begin{cases}
\displaystyle	\frac{1}{\mu (\{v: \langle \nabla  u (x), v \rangle=0 \})} \int_{ \langle \nabla  u (x), v \rangle =0} \langle   D^2  u (x) v, v 
\rangle  d\mu (v) = - C,& \;\; x\in
	  \Omega_0, \\
	u(x) =0,&\;\;x\in\partial \Omega_0.
	\end{cases}
   \end{equation}
   Here
   $$
C =  \frac{1}{2\mu(\{\langle e_N, v \rangle=0\})}\int_{\langle e_N, v \rangle=0} (v_1)^2 \, d \mu (v).
$$
For completeness, we state the precise definition of being a viscosity solution
to \eqref{PDEPb-2} (we remark that \eqref{PDEPb-2-E} and \eqref{PDEPb-2}
are equivalent problems).

\begin{definition} \label{def-2}
An uppersemicontinuous function $\underline{u}$ is a subsolution to \eqref{PDEPb-2}
provided that if
$\underline{u} - \phi$
 has a local maximum at a point $x_0 \in \Omega_0$
for some $\phi \in C^2$, then
$$
\frac{1}{\mu (\{v: \langle \nabla  u (x), v \rangle=0 \})} \int_{ \langle \nabla  u (x), v \rangle =0} \frac12 \langle   D^2  \phi (x_0) v, v 
\rangle  d\mu (v) \geq - C
$$
when $\nabla \phi (x_0) \neq 0$; and 
$$
\Delta \phi (x_0) - \langle  D^2  \phi (x_0) \eta, \eta 
\rangle  \rangle  \geq -1,
$$
for some vector $\eta$ with $|\eta|\leq 1$;
and $\underline{u} \leq 0$ on $\partial \Omega$.

A lowersemicontinuous function $\overline{u}$ is a supersolution to \eqref{PDEPb-2}
provided that if
$\overline{u} - \phi$
 has a local minimum at a point $x_0 \in \Omega_0$
for some $\phi \in C^2$, then
$$
\frac{1}{\mu (\{v: \langle \nabla  u (x), v \rangle=0 \})} \int_{ \langle \nabla  u (x), v \rangle =0} \frac12 \langle   D^2  \phi (x_0) v, v 
\rangle  d\mu (v) \leq - C
$$
when $\nabla \phi (x_0) \neq 0$; and 
$$
\Delta \phi (x_0) - \langle  D^2  \phi (x_0) \eta, \eta 
\rangle  \rangle  \leq -1,
$$
for some vector $\eta$ with $|\eta|\leq 1$;
and $\overline{u} \geq 0$ on $\partial \Omega$.

Finally, a function $u$ is a solution to \eqref{PDEPb-2} if it is both a
super and a subsolution.
\end{definition}

Just for completeness we include a proof of the fact that both definitions are equivalent.

\begin{lemma}
According to Definition~\ref{def-1}, a function $u$ is a viscosity supersolution (subsolution) if and only if it is a viscosity supersolution (subsolution) according to Definition~\ref{def-2}.
\end{lemma}

\begin{proof} The proof is based on the fact that for a smooth function $\phi$
with $\nabla \phi (x) \neq 0$ it holds that (here we are assuming for simplicity that $ \nabla \phi (x)$  points
in the direction of $e_N$),
$$
\begin{array}{l}
\displaystyle
 \frac{1}{\mu(\langle \nabla \phi (x), v \rangle=0)}
\int_{\langle \nabla \phi (x), v \rangle=0} \frac12
\langle D^2 \phi (x) v,v \rangle 
d \mu (v) \\[10pt]
\displaystyle
\qquad =\sum_{i,j=1}^{N-1} \frac{1}{\mu(\langle \nabla \phi (x), v \rangle=0)}
\int_{\langle \nabla \phi (x), v \rangle=0} \frac12
  \frac{\partial^2\phi}{\partial v_i\partial v_j}(x)  v_i,v_j   d \mu (v) \\[10pt]
\displaystyle
\qquad  = C \sum_{i=1}^{N-1}
\frac{\partial^2 \phi}{\partial v_i^2}(x)  = C
\Delta \phi (x) \mid_{\langle \nabla \phi (x), v \rangle=0}
= C \left( \Delta \phi (x) - \Big\langle   D^2  \phi (x) \frac{\nabla \phi}{|\nabla \phi|} (x), 
 \frac{\nabla \phi}{|\nabla \phi|} (x)
\Big\rangle \right)
\end{array}
$$
with 
$$
C= \frac{1}{2\mu (\{ \langle \nabla  u (x), v \rangle=0 \})} \int_{\langle e_N, v \rangle=0} (v_1)^2 \, d \mu (v).
$$

Now, assume that $u$ is a viscosity supersolution according to Definition
\ref{def-1} and that $\phi$ is a smooth function such that 
${u} - \phi$
 has a local minimum at a point $x_0 \in \Omega_0$. 
 Further assume that $\nabla \phi (x_0) \neq 0$ (otherwise there is nothing to
 prove). Then, we have
$$
\Delta \phi (x_0) - \Big\langle   D^2  \phi (x_0) \frac{\nabla \phi}{|\nabla \phi|} (x_0), 
 \frac{\nabla \phi}{|\nabla \phi|} (x_0)
\Big\rangle \leq - 1.
$$
From our previous computation this means that, 
$$
\begin{array}{l}
\displaystyle
 \frac{1}{\mu(\{\langle \nabla \phi (x_0), v \rangle=0\})}
\int_{\langle \nabla \phi (x_0), v \rangle=0} \frac12
\langle D^2 \phi (x_0) v,v \rangle 
d \mu (v) \leq - C .
\end{array}
$$

Conversely, if for a touching test function with non-vanishing gradient we have 
$$
\begin{array}{l}
\displaystyle
 \frac{1}{\mu(\{\langle \nabla \phi (x_0), v \rangle=0\})}
\int_{\langle \nabla \phi (x_0), v \rangle=0} \frac12
\langle D^2 \phi (x_0) v,v \rangle 
d \mu (v) \leq - C,
\end{array}
$$
then, using again our previous computation, we get 
$$
\Delta \phi (x_0) - \Big\langle   D^2  \phi (x_0) \frac{\nabla \phi}{|\nabla \phi|} (x_0), 
 \frac{\nabla \phi}{|\nabla \phi|} (x_0)
\Big\rangle \leq - 1.
$$

The fact that $u$ is a viscosity subsolution to the problem according to Definition
\ref{def-1}
 if and only if it is a viscosity subsolution according to Definition
\ref{def-2} is completely analogous. It can be obtained just reversing the inequalities when appropriate. 
\end{proof}

Now, we collect some well known results for the elliptic PDE problem
\eqref{PDEPb-2}. 
First, a comparison principle holds for viscosity sub and supersolutions to 
 \eqref{PDEPb} and as a consequence we have 
 uniqueness of solutions. Moreover, combining Perron's method with the comparison principle, we obtain also existence of solutions. For the proofs we refer to \cite{ES}.

\begin{theorem}[Comparison Principle]\label{comp.pp}
	If $\underline{u}$ is a  viscosity subsolution and $\overline{u}$ a supersolution to problem
	 \eqref{PDEPb}, then $$\underline{u} (x)\leq \overline{u} (x) \qquad \mbox{ in }\overline{\Omega}_0.$$
{\rm (Existence and uniqueness for the PDE)} There exists a unique solution $u$ to   \eqref{PDEPb}.  Moreover, $u$ is a continuous function in $\overline{\Omega}_0$.
\end{theorem} 

Let us point out that when the domain is a large ball $\Omega_0 = B_R (0)$ the solution to 
	\begin{equation}
		\label{PDEPb-2-E-99}
\begin{cases}
\displaystyle	\Delta u(x) - \Big\langle   D^2  u (x) \frac{\nabla u}{|\nabla u|} (x), 
 \frac{\nabla u}{|\nabla u|} (x)
\Big \rangle  = - L,& \;\; x\in
	  \Omega_0, \\
	u(x) =0,&\;\;x\in\partial \Omega_0,
	\end{cases}
   \end{equation}
   is explicit. Computing the Laplacian in \eqref{PDEPb-2-E-99} in radial coordinates we find that the solution is given by
   $$
   u(x) = \frac{L}{2(N-1)} (R^2-|x|^2).
   $$
   Hence, taking $L>1$ and $R$ large such that $\Omega_0 \subset\subset B_R(0)$ we get a strict supersolution to our problem \eqref{PDEPb-2-E}.

We will also need a very useful tool from probability theory.  

\subsection{Probability. The Optional Stopping Theorem.}
We briefly recall (see \cite{Williams}) that a sequence of random variables
$\{M_{k}\}_{k\geq 1}$ is a supermartingale (submartingales) if
$$ \E[M_{k+1}\arrowvert M_{0},M_{1},...,M_{k}]\leq M_{k} \ \ (\geq).$$
Then, the Optional Stopping Theorem, that we will call {\it (OSTh)} in what follows, says:
given $\tau$ a stopping time such that one of the following conditions hold,
\begin{itemize}
\item[(a)] The stopping time $\tau$ is bounded almost surely;
\item[(b)] It holds that $\E[\tau]<\infty$ and there exists a constant $c>0$ such that $$\E[M_{k+1}-M_{k}\arrowvert M_{0},...,M_{k}]\leq c;$$
\item[(c)] There exists a constant $c>0$ such that $|M_{\min \{\tau,k\}}|\leq c$ almost surely for every $k$.
\end{itemize}
Then, 
$$ \E[M_{\tau}]\leq \E [M_{0}] \ \ (\geq)$$
if $\{M_{k}\}_{k\geq 0}$ is a supermartingale (submartingale).

For the proof of this classical result we refer to \cite{Doob,Williams}.

\subsection{A technical lemma} Now, we prove a technical lemma that will be the key ingredient to pass
to the limit in the DPP and obtain viscosity solutions to the mean curvature
equation. In the next statement recall that $\sigma$ and $\mu$ denote the surface measure on
the spheres $S^{N-1}$ and $S^{N-2}$ respectively. Here, we use $v=(\tilde{v},v_N) \in {S}^{N-1}$ to denote a unitary vector with $\tilde{v}\in \mathbb{R}^{N-1}$ and $v_N:=\langle v, e_N\rangle$, the last coordinate of $v$.

\begin{lemma} \label{lema.clave}  Let $\delta_\eps = \eps^{1/2}$. Take 
two families of sets, $\{B_\eps^\ast\}_{\eps>0}\subset S^{N-1}$ such that
$$
	B_\eps^\ast=\Big\{v\in S^{N-1} \ : \ \langle e_N,v\rangle \le \theta_\eps\Big\},
	$$
	with $\theta_\eps$ such that
$\sigma (B_\eps^\ast) = \frac12\sigma (S^{N-1})+\delta_\eps$,
and $\{A_\eps^\ast\}_{\eps>0}\subset S^{N-1}$ any family verifying $\sigma (A_\eps^\ast)\ge \frac12\sigma (S^{N-1})+\delta_\eps$.
Assume that
	$$
	-\eps\le\frac{1}{\sigma (A_\eps^\ast\cap B_\eps^\ast)}\int_{A_\eps^\ast\cap B_\eps^\ast}\langle e_N,v\rangle d\sigma (v)\le 0 .
	$$
	Then,
	given a continuous function $f:S^{N-1}\rightarrow\R$, it holds that
	\begin{equation}
		\displaystyle \frac{1}{\sigma (A_\eps^\ast\cap B_\eps^\ast)}\int_{A_\eps^\ast\cap B_\eps^\ast}f(v) d\sigma (v)\longrightarrow \frac{1}{\mu (S^{N-2})}\int_{S^{N-1}\cap\{\langle e_N,v\rangle=0\}}f((\widetilde{v},0)) d\mu (\widetilde{v})
	\end{equation}
	when $\eps\to 0$. 
\end{lemma}

\begin{proof} 
First, we notice that since $\delta_\eps \approx \eps^{1/2}$
we have $\theta_\eps \approx \eps^{1/2}$.

Let us define the set $U_\eps=\{v\in S^{N-1} \ : \ |\langle e_N,v\rangle|\le \theta_\eps\}$ and consider
	$$
	\begin{array}{ll}
		\displaystyle \left| \vint_{A_\eps^\ast\cap B_\eps^\ast}f(v) d\sigma (v)-\vint_{S^{N-1}\cap\{\langle e_N,v\rangle=0\}}f((\widetilde{v},0)) d\mu (\widetilde{v})\right| \\[12pt]
		\displaystyle \le \underbrace{\left| \vint_{A_\eps^\ast\cap B_\eps^\ast}f(v) d\sigma (v)-\vint_{U_\eps}f(v) d\sigma (v)\right|}_{I}\\
		\displaystyle \qquad +\underbrace{\left| \vint_{U_\eps}f(v) d\sigma(v)-\vint_{S^{N-1}\cap\{\langle e_N,v\rangle=0\}}f((\widetilde{v},0)) d\mu (\widetilde{v})\right|}_{II}.
	\end{array}
	$$
	
	Let us estimate the second term, $II$. It holds that
	$$
	\lim_{\eps\to 0}\frac{\sigma (U_\eps)}{2\theta_\eps} = \mu (S^{N-2}).
	$$
	Using that $f$ is continuous, by Lebesgue differentiation theorem, we have 
	$$
	\frac{1}{2\theta_\eps}\int_{\{-\theta_\eps\le v_N\le\theta_\eps\}}f((\widetilde{v},v_N))d{v_N}\longrightarrow f((\widetilde{v},0))
	$$
	when $\eps\to 0$. Hence, we obtain 
	$$
	\begin{array}{l}
	\displaystyle 
	\lim_{\eps \to 0}
	\frac{2\theta_\eps}{\sigma (U_\eps)} \frac{1}{2\theta_\eps} \int_{U_\eps}
	f({v})d\sigma (v) \\[10pt]
	\qquad \displaystyle = \frac{1}{\mu (S^{N-2})}\int_{S^{N-2}}f((\widetilde{v},0))d\mu ((\widetilde{v})).
\end{array}	
	$$
	Thus, $II\to 0$ when $\eps\to 0$.
	
	Now, let us consider $I$, using that $\lVert f\rVert_\infty \le C$, we get
	$$
	\begin{array}{ll}
		\displaystyle \left| \vint_{A_\eps^\ast\cap B_\eps^\ast}f(v) d\sigma (v)-\vint_{U_\eps}f(v) d\sigma (v)\right|\\[10pt]
		\displaystyle \qquad	\le \frac{C}{\sigma (U_\eps)}\Big[\sigma (A_\eps^\ast\cap\{\langle e_N,v\rangle<-\theta_\eps\})+\sigma ((A_\eps^\ast)^c\cap U_\eps)\Big]. 
	\end{array}
	$$
	Here we used that $\sigma(A_\eps^\ast\cap B_\eps^\ast)\geq \sigma (U_\eps)$. We will prove that 
	$$
	\frac{\sigma (A_\eps^\ast\cap\{\langle e_N,v\rangle<-\theta_\eps\})}{\sigma (U_\eps)}\rightarrow 0
	$$
	when $\eps \to 0$. Suppose that this is not true, then there exist a positive constant $L$ and a sequence $\eps_j\rightarrow 0$ such that 
	$$
	\frac{\sigma (A_{\eps_j}^\ast\cap\{\langle e_N,v\rangle<-\theta_{\eps_j}\})}{\sigma (U_{\eps_j})}\ge L> 0
	$$
	for all $j\in\mathbb{N}$.
	Then 
	$$
	\frac{\sigma (A_{\eps_j}^\ast\cap\{\langle e_N,v\rangle<-(1+\alpha)\theta_{\eps_j}\})}{\sigma (U_{\eps_j})}\ge \frac{L}{2}> 0,
	$$
	for some $\alpha>0$ and for all $j\geq j_0$ for some $j_0\in\mathbb{N}$. From this point, we omit the subindex $j$
		for clarity. We can assume that $\sigma (A_\eps^\ast \cap \{-(1+\alpha)\theta_\eps<\langle e_N,v\rangle<-\theta_\eps\})=0$. Let us use the hypothesis to obtain
	$$
	\begin{array}{ll}
		\displaystyle -\eps\le\frac{1}{\sigma (A_\eps^\ast\cap B_\eps^\ast)}\int_{A_\eps^\ast\cap B_\eps^\ast}\langle e_N,v\rangle d\sigma (v) \leq \frac{1}{\sigma (U_\eps)}\int_{A_\eps^\ast\cap B_\eps^\ast}\langle e_N,v\rangle d\sigma (v) \\[12pt]
		\displaystyle \qquad \le \frac{1}{\sigma (U_\eps)}\Big[\int_{A_\eps^\ast\cap \{-\theta_\eps\le \langle e_N,v\rangle\le \theta_\eps\}}\langle e_N,v\rangle d\sigma (v)+\int_{A_\eps^\ast\cap \{ \langle e_N,v\rangle<- (1+\alpha)\theta_\eps\}}\langle e_N,v\rangle d\sigma (v)\Big]\\[12pt]
		\displaystyle \qquad \le \delta_\eps\frac{\sigma ( A_\eps^\ast\cap\{\langle e_N,v\rangle<-(1+\alpha)\theta_\eps\})}{\sigma (U_\eps)}-(1+\alpha)\theta_\eps\frac{\sigma ( A_\eps^\ast\cap\{\langle e_N,v\rangle<-(1+\alpha)\theta_\eps\})}{\sigma (U_\eps)}\\[12pt]
		\displaystyle \qquad \le-\alpha\theta_\eps\frac{\sigma ( A_\eps^\ast\cap\{\langle e_N,v\rangle<-(1+\alpha)\theta_\eps\})}{\sigma (U_\eps)} \le-\alpha\theta_\eps\frac{L}{2}
	\end{array}
	$$
	which is a contradiction since 
	we have
	$$1 \leq - \alpha \frac{L}{2} \frac{\theta_\eps}{\eps}\to -\infty, \qquad \mbox{as } \eps \to 0,$$ 
	because $\delta_\eps=\eps^{1/2}$ and thus $\theta_\eps \approx \eps^{1/2}$. In the previous estimates we used 
	that
	\[
	0=\int_{U_\eps}\langle e_N,v\rangle d\sigma (v) =\int_{A_\eps^\ast\cap \{-\theta_\eps\le \langle e_N,v\rangle\le \theta_\eps\}}\langle e_N,v\rangle d\sigma (v)+\int_{(A_\eps^\ast)^c \cap \{-\theta_\eps\le \langle e_N,v\rangle\le \theta_\eps\}}\langle e_N,v\rangle d\sigma (v)
	\]
	implies
		\[
	\begin{array}{l}
		\displaystyle 	\frac{1}{\sigma (A_\eps^\ast\cap B_\eps^\ast)}\int_{A_\eps^\ast\cap \{-\theta_\eps\le \langle e_N,v\rangle\le \theta_\eps\}}\langle e_N,v\rangle d\sigma (v)\\[12pt]
		\qquad \displaystyle  =\frac{1}{\sigma (A_\eps^\ast\cap B_\eps^\ast)}\int_{(A_\eps^\ast)^c \cap \{-\theta_\eps\le \langle e_N,v\rangle\le \theta_\eps\}}-\langle e_N,v\rangle d\sigma (v)\\[12pt]
		\displaystyle \qquad \le\theta_\eps\frac{\sigma ((A_\eps^\ast)^c \cap \{-\theta_\eps\le \langle e_N,v\rangle
		  \le \theta_\eps\})}{\sigma (U_\eps)}
		  \\[12pt]
		\qquad \displaystyle =\theta_\eps\frac{\sigma ( A_\eps^\ast\cap\{\langle e_N,v\rangle<-(1+\alpha)\theta_\eps\})}{\sigma (U_\eps)}.
	\end{array}
	\]
	
	Thus, we also have that $I\to 0$ when $\eps\to 0$. To show that
	$$
	\frac{\sigma ((A_\eps^\ast)^c\cap U_\eps)}{\sigma (U_\eps)}\rightarrow 0
	$$
	is analogous. This ends the proof.
\end{proof}

\section{The value of the game is the unique solution to the DPP}\label{sect.Game}

In this section we show that the value of the game is characterized as the unique solution to the DPP. This fact is crucial when we pass to the limit in the viscosity sense and obtain the limit is a solution to the mean curvature equation. 

\begin{lemma} \label{lema-DPP} There exists a solution to
\begin{equation}
	\label{DPP-22} 
		\left\{
	\begin{array}{ll}
	\displaystyle
	v(x) = \sup_{A} \inf_{B}  \left\{ \frac{1}{\sigma(A\cap B)}
\int_{A\cap B}  v \big(x + v\varepsilon \big) 
d \sigma (v)\right\} + \varepsilon^2K,  & x\in \Omega_0,\\[10pt]
	v(x) = 0 ,  & x\in \mathbb{R}^N \setminus 
\Omega_0.
	\end{array} \right.
\end{equation}
\end{lemma}

\begin{proof}
We start with $w_0 \equiv 0$, and we define inductively
$$
w_{n+1}(x) = \sup_{A} \inf_{B}  \left\{ \frac{1}{\sigma(A\cap B)}
\int_{A\cap B}  w_n \big(x + v\varepsilon \big) 
d \sigma (v)\right\} +  \varepsilon^2K,  \qquad  \mbox{for } x\in \Omega_0
$$
with $w_{n+1} \equiv 0$ for $x\in \mathbb{R}^N \setminus 
\Omega_0.$ 

First, note that for each $x$ in the interior of $\Omega_0$, the chosen sets $A$ and $B$ where supremum and infimum are attained are not necessarily contained in $\Omega_0$. In fact, for points close enough to the boundary, we expect that a most part of the selected $B$ lies outsides $\Omega_0$. Consequently, $w_n$ is not constant within $\Omega_0.$

Second, it is not hard to check by an inductive argument that $w_n$ is increasing with $n$, and thus we can consider
$$
v(x) = \lim_{n\to \infty} w_n (x).
$$
Notice that 
$$
   u(x) = \frac{L}{2(N-1)} (R^2-|x|^2).
   $$
  with $L>1$ and $R$ large such that $\Omega_0 \subset\subset B_R(0)$ is a strict supersolution to our problem \eqref{PDEPb-2-E}. Therefore, by a simple Taylor expansion like the one that we outline in the 
  introduction we have that $u(x)$ is a strict supersolution to \eqref{DPP-22}. This argument provides
  a uniform upper bound for the sequence $w_n$.

Now, we observe that,  by the monotone convergence theorem we have 
$$
\begin{array}{l}
\displaystyle
\lim_{n\to \infty} \sup_{A} \inf_{B}  \left\{ \frac{1}{\sigma(A\cap B)}
\int_{A\cap B}  w_n \big(x + v\varepsilon \big) 
d \sigma (v)\right\} + \varepsilon^2 K \\[10pt]
\displaystyle =
 \sup_{A} \inf_{B}  \left\{ \frac{1}{\sigma(A\cap B)}
\int_{A\cap B}  v \big(x + v\varepsilon \big) 
d \sigma (v)\right\} +  \varepsilon^2 K.
\end{array}
$$
Hence, we conclude that
$$
v(x) = \sup_{A} \inf_{B}  \left\{ \frac{1}{\sigma(A\cap B)}
\int_{A\cap B} v \big(x + v\varepsilon \big) 
d \sigma (v)\right\} +  \varepsilon^2K,  \qquad  \mbox{for } x\in \Omega_0
$$
with $v \equiv 0$ for $x\in \mathbb{R}^N \setminus 
\Omega_0$. Thus, we have the desired existence of a solution to the DPP.
\end{proof}

Recall that the value of the game for Paul is given by 
\[\overline{u}^\varepsilon (x_0) \! = \inf_{S_c} \sup_{S_p} \mathbb{E}^{x_0}_{S_p,S_c} 
\Big[ \varepsilon^2 K \times (\mbox{number of plays}) \Big].
\]
This is the best value that Carol may guarantee to obtain.

Analogously, the value of the game from Carol's viewpoint is
\[\underline{u}^\varepsilon (x_0) \! =  \sup_{S_p} \inf_{S_c} \mathbb{E}^{x_0}_{S_p,S_c} 
\Big[ \varepsilon^2 K \times (\mbox{number of plays}) \Big].
\]
We say that the game has a value if we can reverse $\inf$ with $\sup$, that is, when 
\[\overline{u}^\varepsilon (x_0) = \underline{u}^\varepsilon (x_0) := u^\varepsilon (x_0). \]
Our next result shows that this game has a value and, moreover, proves that
the value function is the solution to the DPP.

\begin{theorem} \label{teo-valor-sol-DPP}
The game has a value that is characterized as the unique solution to
\eqref{DPP-22}.
\end{theorem}

\begin{proof}
In Lemma \ref{lema-DPP} we proved that the DPP has a solution $v(x)$. In what follows we use this solution to build quasi-optimal strategies for the players and use probabilistic arguments to show that indeed the value of the game 
coincides with the solution to the DPP, $v$.

At every position of the game, $x_k$, Paul (who wants to maximize) chooses $A^*$ in such a way that
$$
\begin{array}{l}
\displaystyle
\sup_{A} \inf_{B}  \left\{ \frac{1}{\sigma(A\cap B)}
\int_{A\cap B}  v \big(x_k + v\varepsilon \big) 
d \sigma (v)\right\} \\[10pt]
\qquad \displaystyle \leq  \inf_{B}  \left\{ \frac{1}{\sigma(A^*\cap B)}
\int_{A^*\cap B}  v \big(x_k + v\varepsilon \big) 
d \sigma (v)\right\} + \frac{\eta}{2^{k+1}}.
\end{array}
$$

Given this strategy for Paul and any strategy $S_{c}$ for Carol we consider the sequence of random variables
given by
$$
M_{k}=
 \displaystyle v(x_{k})-\frac{\eta}{2^{k}} 
$$
Let us see that $(M_{k})_{k\geq 0}$ is a submartingale. 
To this end we need to estimate 
$$
\E_{S_{p}^{*},S_{c}}[M_{k+1}\mid_{M_k}].
$$

Since we are using the strategies $S_{p}^{*}$ and $S_{c}$, it holds that 
\begin{align*}
    \E_{S_{p}^{*},S_{c}}[M_{k+1}\mid_{M_k}]
    &= \E_{S_{p}^{*},S_{c}}[ \displaystyle v(x_{k+1})-\frac{\eta}{2^{k+1}} \mid_{M_k}] \\
    &= \left\{ \frac{1}{\sigma(A^*\cap B)} \int_{A^*\cap B}  v \big(x_k + v\varepsilon \big) d \sigma (v)\right\} - \frac{\eta}{2^{k+1}}.
\end{align*}
Here $A^*$ corresponds to Paul's choice and $B$ to Carol's. 
Now, we have
\[
\begin{array}{l}
\displaystyle 
\sup_{A} \inf_{B}  \left\{ \frac{1}{\sigma(A\cap B)}
\int_{A\cap B}  v \big(x_k + v\varepsilon \big) 
d \sigma (v)\right\} \\[10pt] 
\displaystyle \qquad \leq  \inf_{B}  \left\{ \frac{1}{\sigma(A^*\cap B)}
\int_{A^*\cap B}  v \big(x_k + v\varepsilon \big) 
d \sigma (v)\right\} + \frac{\eta}{2^{k+1}} \\[10pt] 
\displaystyle \qquad \leq  \left\{ \frac{1}{\sigma(A^*\cap B)}
\int_{A^*\cap B}  v \big(x_k + v\varepsilon \big) 
d \sigma (v)\right\}  + \frac{\eta}{2^{k+1}}.
\end{array}
\]

Therefore, we arrive to
$$
\E_{S_{p}^{*},S_{c}}[M_{k+1}\mid_{M_k}]
\geq \sup_{A} \inf_{B}  \left\{ \frac{1}{\sigma(A\cap B)}
\int_{A\cap B}  v \big(x_k + v\varepsilon \big) 
d \sigma (v)\right\}  - \frac{\eta}{2^k}.
$$

As $v$ is a solution to the DPP \eqref{DPP} we obtain 
$$
\E_{S_{p}^{*},S_{c}}[M_{k+1}\mid_{M_k}]\geq v(x_{k})-\frac{\eta}{2^{k}}=M_{k}
$$
as we wanted to show.

Therefore, $(M_{k})_{k\geq 0}$ is a submartingale. Using the optional stopping theorem 
(recall that we have that $\tau$ is finite a.s. and that we have that $M_k$ is uniformly bounded we conclude that 
$$
\E_{S_{p}^{*},S_{c}}[M_{\tau}]\geq M_{0}
$$
where $\tau$ is the first time such that $x_{\tau}\notin\Omega$. Then, 
$$
\E_{S_{p}^{*},S_{c}}[\eps^2K\tau]\geq v (x_{0})-\eta.
$$
We can compute the infimum  in $S_c$ and then
the supremum in $S_p$ to obtain
$$
\sup_{S_p}\inf_{S_c}\E_{S_p,S_c}[\eps^2 K \tau]\geq v(x_{0})-\eta.
$$

An analogous computation shows that
$$
\inf_{S_{c}}\sup_{S_{p}}\E_{S_{p},S_{c}}[\eps^2 K\tau]\leq v(x_{0})+\eta.
$$

To end the proof we just observe that 
$$
\sup_{S_{p}}\inf_{S_{c}}\E_{S_{p},S_{c}}[\eps^2 K\tau]\leq \inf_{S_{c}}\sup_{S_{p}}\E_{S_{p},S_{c}}[\eps^2K\tau].
$$
Therefore,
$$
v(x_{0})-\eta\leq\sup_{S_{p}}\inf_{S_{c}}\E_{S_{p},S_{c}}[\eps^2 K\tau]\leq \inf_{S_{c}}\sup_{S_{p}}\E_{S_{p},S_{c}}[\eps^2K\tau]\leq v(x_{0})+\eta
$$
Since $\eta>0$ is arbitrary, we conclude that the game has a value 
and that 
$$
v(x_{0})=\sup_{S_{p}}\inf_{S_{c}}\E_{S_{p},S_{c}}[\eps^2 K\tau]= \inf_{S_{c}}\sup_{S_{p}}\E_{S_{p},S_{c}}[\eps^2K\tau].
$$
That is
\[
v(x_0)=\underline{u}^\eps(x_0)=\overline{u}^\eps(x_0).
\]
This ends the proof. 
\end{proof}

As a consequence of the previous arguments we obtain a comparison principle for sub and super solutions to the DPP. First, let us state what we understand by a super and a subsolution to the DPP. 

\begin{definition}
A function $\overline{u}^\varepsilon$ is a supersolution to the DPP if it verifies
\begin{equation}
	\label{DPP-22-super} 
		\left\{
	\begin{array}{ll}
	\displaystyle
	\overline{u}^\varepsilon (x) \geq \sup_{A} \inf_{B}  \left\{ \frac{1}{\sigma(A\cap B)}
\int_{A\cap B}  \overline{u}^\varepsilon \big(x + v\varepsilon \big) 
d \sigma (v)\right\} +  \varepsilon^2K,  & x\in \Omega_0,\\[10pt]
	\overline{u}^\varepsilon (x) \geq 0 ,  & x\in \mathbb{R}^N \setminus 
\Omega_0.
	\end{array} \right.
\end{equation}

A function $\underline{u}^\varepsilon$ is a subsolution to the DPP 
when the reverse inequalities hold, that is, when 
\begin{equation}
	\label{DPP-22-sub} 
		\left\{
	\begin{array}{ll}
	\displaystyle
	\underline{u}^\varepsilon (x) \leq \sup_{A} \inf_{B}  \left\{ \frac{1}{\sigma(A\cap B)}
\int_{A\cap B}  \underline{u}^\varepsilon \big(x + v\varepsilon \big) 
d \sigma (v)\right\} +  \varepsilon^2K,  & x\in \Omega_0,\\[10pt]
	\underline{u}^\varepsilon (x) \leq 0 ,  & x\in \mathbb{R}^N \setminus 
\Omega_0.
	\end{array} \right.
\end{equation}
\end{definition}

Now, we are ready to prove the comparison argument.

\begin{theorem} \label{compar-DPP}
Let $\overline{u}^\varepsilon$ be a supersolution to the DPP and $\underline{u}$ a subsolution. Then
$$
\overline{u}^\varepsilon (x) \geq \underline{u}^\varepsilon  (x), \qquad \mbox{ for every } 
x \in \mathbb{R}^N.
$$
\end{theorem}

\begin{proof} We only have to observe that the previous proof shows that
$\overline{u}^\varepsilon (x)$ is bigger or equal than the value of the game.
Indeed, one can follow the same argument as before to show that
$$
M_{k}=
 \displaystyle \overline{u}^{\varepsilon}(x_{k})+\frac{\eta}{2^{k}} 
$$
is a supermartingale
using that $ \overline{u}^{\varepsilon}$ verifies \eqref{DPP-22-super}.

Analogously, one can also show that $\underline{u}$ is less or equal than the value of the game. 

Then, we conclude that 
$$
\overline{u}^\varepsilon (x) \geq \underline{u}^\varepsilon  (x), \qquad \mbox{ for every } 
x \in \mathbb{R}^N
$$
as we wanted to show.
\end{proof}

\section{The limit of the value functions for the game is the solution
to the mean curvature equation}\label{sect.GameAndEq}

In this section, we will obtain that the unique solution $u$ to \eqref{PDEPb} is the locally uniform limit of the value functions $u^{\varepsilon}$ of the $\varepsilon$-game as $\varepsilon$ goes to zero, see Theorem  \ref{maintheorem.juego.conv}.
Our strategy is to prove that the functions 
\begin{equation} \label{subsuperviscositysols}
	\overline{u}(x) = \limsup_{\substack{\varepsilon \to 0^+ \\ y \to x}} u^\varepsilon(y), \quad \underline{u}(x) = \liminf_{\substack{\varepsilon \to 0^+ \\ y \to x}} u^\varepsilon(y),
\end{equation} 
are,  respectively, a viscosity subsolution and a viscosity supersolution to the problem \eqref{PDEPb} with $\overline{u}(x) = \underline{u}(x) =0$ for $x\in \partial \Omega_0$. Therefore, by the Comparison Principle, Theorem \ref{comp.pp}, we get that $\overline{u}\leq\underline{u}.$ Since by definition $\underline{u}\leq\overline{u}$, we then obtain that $\overline{u}=\underline{u}$. Hence, the limit of the family $u^{\varepsilon}$ exists, and it is the unique solution $u$ to \eqref{PDEPb}. Finally, as a consequence of the convergence
$u^{\varepsilon}\to u$, for each $t>0$ we will be able to approximate the set where $u(\cdot)$ is 
bigger than $t$ by the set where $u^{\varepsilon}(\cdot)$ 
is bigger than $t$ for $\varepsilon$ small enough, see Corollary~\ref{cor1}. 

\subsection{Uniform bounds} First, we show that the value function of the game, $u^\eps$ is uniformly bounded. 

\begin{lemma} \label{lema.cota.uniforme}
There exists a constant $C$, independent of $\eps$, such that
\begin{equation} \label{cota.uniforme}
\sup_{x\in \Omega_0 } |u^\eps (x) |\leq C.
\end{equation}
\end{lemma}

\begin{proof}
We have that $u^\eps (x) \geq 0$. Hence, we need to show only an upper bound. 
To simplify the notation we assume that this point is the origin $z=0$. Assume that at any point $x\in \Omega_0$ Carol chooses $B_\eps$ as the set 
$$B_{\varepsilon}=\Big\{  v\in S^{N-1}: \langle v, \frac{x}{|x|}  \rangle \geq  - \theta_\varepsilon\Big\}$$
with $\theta_\eps$ such that 
$\sigma (B_\eps) = \frac12\sigma (S^{N-1})+\delta_\eps$.
That is, Carol aims to go as far as possible from $z=0$. 

Then, for any set $A$ that Paul may choose, we have that the next position of the game
provided that we arrive to $x_k$ satisfies
$$
\begin{array}{ll}
\displaystyle \mathbb{E} [|x_{k+1}|^2\mid_{x_k}]=\frac{1}{\sigma(A\cap B_\eps)}
\int_{A\cap B_\eps}  |x_k + v\varepsilon |^2
d \sigma (v) \\
\displaystyle = \frac{1}{\sigma(A\cap B_\eps)}
\int_{A\cap \{-\theta_\eps\le\langle v,\frac{x_k}{|x_k|}\rangle\le\theta_\eps \}}  |x_k + v\varepsilon |^2
d \sigma (v) \\
\displaystyle \qquad +\frac{1}{\sigma(A\cap B_\eps)}
\int_{A\cap \{\theta_\eps\le\langle v,\frac{x_k}{|x_k|}\rangle\}}  |x_k + v\varepsilon |^2
d \sigma (v)\\
\displaystyle \ge\frac{1}{\sigma(A\cap B_\eps)}
\int_{A\cap \{-\theta_\eps\le\langle v,\frac{x_k}{|x_k|}\rangle\le\theta_\eps \}}  |x_k + v\varepsilon |^2
d \sigma (v) \\
\displaystyle \qquad +\frac{1}{\sigma(A\cap B_\eps)}
\int_{A^c\cap \{-\theta_\eps\le\langle v,\frac{x_k}{|x_k|}\rangle\le\theta_\eps \}}  |x_k + v\varepsilon |^2
d \sigma (v)\\
\displaystyle = \frac{1}{\sigma(A\cap B_\eps)}
\int_{\{-\theta_\eps\le\langle v,\frac{x_k}{|x_k|}\rangle\le\theta_\eps \}}  |x_k + v\varepsilon |^2
d \sigma (v)= |x_k|^2 +  \varepsilon^2.
\end{array}
$$
Here we used that 
\[
\sigma(A\cap \{\theta_\eps\le\langle v,\frac{x_k}{|x_k|}\rangle\})\ge\sigma( A^c\cap \{-\theta_\eps\le\langle v,\frac{x_k}{|x_k|}\rangle\le\theta_\eps \})
\]
and 
\[
 |x_k + v\varepsilon |^2\mid_{ A\cap \{\theta_\eps\le\langle v,\frac{x_k}{|x_k|}\rangle\}}\ge |x_k + v\varepsilon |^2\mid_{A^c\cap \{-\theta_\eps\le\langle v,\frac{x_k}{|x_k|}\rangle\le\theta_\eps \}}.
\]
Then
\[
\frac{1}{\sigma(A\cap B_\eps)}
\int_{A\cap \{\theta_\eps\le\langle v,\frac{x_k}{|x_k|}\rangle\}}  |x_k + v\varepsilon |^2
d \sigma (v)\ge \frac{1}{\sigma(A\cap B_\eps)}
\int_{A^c\cap \{-\theta_\eps\le\langle v,\frac{x_k}{|x_k|}\rangle\le\theta_\eps \}}  |x_k + v\varepsilon |^2
d \sigma (v)
\]

Now, if we introduce the sequence of random variables 
$$M_k =   |x_{k}|^2 -  \varepsilon^2 k .$$ 
It holds that
$$
\mathbb{E} [M_{k+1}\mid_{M_k}] \geq 
|x_{k}|^2+\eps^2 -  \varepsilon^2 (k+1) = M_k.
$$
Then, we have that $M_k$ is a submartingale and using the
optional stopping theorem (OSTh), we obtain
$$
\mathbb{E} [M_{\tau}] \geq 
M_0 = |x_{0}|^2
$$
Hence we conclude that
$$
  \mathbb{E} [\eps^2\tau ] \leq C(\Omega)
$$
This says that the expected value for the number of plays 
when Carol uses the previous strategy is bounded by $C/\varepsilon^2$.
Hence, the value of the game verifies
$$
u^\varepsilon (x_0) \leq \sup_{S_p} \mathbb{E}_{S_p,S^*_c} [\varepsilon^2 \tau]
\leq C
$$
as we wanted to show.

As an alternative argument, we can use that 
$$
   u(x) = \frac{L}{2(N-1)} (R^2-|x|^2).
   $$
  with $L>1$ and $R$ large such that $\Omega_0 \subset\subset B_R(0)$ is a strict supersolution to our problem \eqref{PDEPb-2-E}. Then, by a simple Taylor expansion like the one that we outline in the 
  introduction we have that $u(x)$ is a strict supersolution to \eqref{DPP-22} for $\varepsilon $ small
  enough. Then, the comparison principle for the DPP, Theorem \ref{compar-DPP}, provides
  a uniform upper bound for the solution to the DPP $u^\varepsilon$ for $\varepsilon$ small enough.
\end{proof}

\subsection{Estimates near the boundary}
 Next, we prove that the half-relaxed limits 
 of the family $u^\eps$, $\overline{u}$ and $\underline{u}$, given by \eqref{subsuperviscositysols},
 attain the boundary condition.

\begin{proposition} \label{ppDatoInicial}
	Let $\overline{u}$ and $\underline{u}$ be defined as in \eqref{subsuperviscositysols}. Then  $$\overline{u}(x) = \underline{u}(x) =0, \qquad \mbox{ for every } x\in \partial \Omega_0.$$
\end{proposition}

\begin{proof} 
	We will assume, to simplify the notation, that we start the game at a point $x_0$ of the form $x_0=(0,..,0,x_N^0) \in \Omega_0$
	with $x_N^0>0$ that is close to
	$y_0=0 \in \partial \Omega_0$.
	
	Now, let us show the two inequalities needed for the proof. 
	The lower bound is straightforward since we have that
	\begin{equation}
		u_\varepsilon (x_0) \geq 0,
	\end{equation}
	and then we conclude that
	$$
	\overline{u}(x) \geq \underline{u}(x) \geq 0
	$$
	in the whole $\overline{\Omega}_0$. 
	
	To find a reverse estimate is more delicate. To this end we
	need to choose a strategy for Carol (the player who wants 
	to minimize the expected payoff), associated to the following claim.
	
	CLAIM: Given $\eta>0$, there exist $r_0>0$ and $z\in\R^N$ such that if $|x_0|<r_0$, and $x_\tau\in\partial\Omega_0\cap(B_{|x_0-z|}(z))^c$ for some $\tau\in\mathbb{N}$, we get
	\[
	\Big||x_\tau-z|^2-|x_0-z|^2\Big|<\eta.
	\]
	Suppose that the claim is true (we will prove it later). At each turn, if the position of the token is in $x_k$ we let Paul choose freely while Carol chooses $B^\ast$ as the set 
	$$B^\ast=\Big\{  v\in S^{n-1} : \langle v, \frac{x_k-z}{|x_k-z|} \rangle \geq  -\theta_\varepsilon\Big\},$$
	with $\theta_\eps$ such that
	$\sigma (B_\eps^\ast) = \frac12\sigma (S^{N-1})+\delta_\eps$.
	
	Since Carol aims to minimize the expected payoff, we get a bound from above for the value of the game, it holds that
	\begin{equation}
		u_\varepsilon (x_0) \leq \sup_{S_p} \mathbb{E}^{x_0}_{S_c^\ast, S_p} [\varepsilon^2 \tau].
	\end{equation}
	
	Then, for any sets $A$ (any strategy for Paul's choices), let us consider the following sequence of random variables 
	\[
	M_k=|x_k-z|^2.
	\]
	We have
	$$
	\begin{array}{ll}
\displaystyle \mathbb{E} [M_{k+1} \mid_{x_k}]=	\mathbb{E} [|x_{k+1}-z|^2 \mid_{x_k}] =
\left\{ \frac{1}{\sigma_{n-1}(A\cap B^\ast)}
\int_{A\cap B^\ast}  |x_k-z+\eps v|^2
d \sigma_{n-1} (v)\right\} \\[10pt]
       \displaystyle \geq |x_k-z|^2+\eps^2 =M_k+ \eps^2\ge M_k.
	\end{array}
	$$

	Thus, using the optional stopping theorem we get
	\[
	\E[M_\tau]\ge M_0,
	\]
	that is
	\[
	\E[|x_\tau-z|^2]\ge |x_0-z|^2.
	\]
	This implies that $x_\tau\in\partial\Omega_0\cap(B_{|x_0-z|}(z))^c$. On the other hand, if we consider
	\[
	N_k=M_k-k\eps^2
	\]
	we obtain
	\[
	\E[N_{k+1}\mid_{x_k}]=\E[M_{k+1}-(k+1)\eps^2\mid_{x_k}]\geq M_k+\eps^2-(k+1)\eps^2=N_k,
	\]
	that is, $N_k$ is a submartingale. Then, if we use the optional stopping theorem we get
	\[
	\E[N_\tau]\ge N_0.
	\]
	That is
	\[
	\E[|x_\tau-z|^2-\tau\eps^2]\ge |x_0-z|^2.
	\]
	Finally, we obtain
	\[
	\E[\tau\eps^2]\le \E[|x_\tau-z|^2]-|x_0-z|^2
	\]
	which is small, using the claim. This implies that
	 \[
	 \uline{u}(0)=\liminf u_\eps(x)=0 \quad \mbox{and} \quad \overline{u}(0)=\limsup u_\eps(x)=0.
	 \]
	 
	 Let us prove the claim.
	 
	 \textit{Proof of the claim.} We can assume that $\overline{\Omega}_0\subset U$ where
	 \[
	 U=\Big\{x\in\R^N : x_N\ge c\sum_{j=1}^{N-1}x_j^2 \Big\},
	 \]
	 with $c>0$. Here we use the fact that we assumed that $\Omega_0$ is strictly convex. We can also suppose that $\partial\Omega_0\cap \partial U=\{0\}$. We will work with the set $U$. If we consider $z=(0,\dots,\Lambda)$ with $\Lambda>0$ large enough. One can check that the function
	 \[
	 \mu: U\cap\{x: \langle x-x_0,e_N\rangle\le 0 \},
	 \]
	 defined by $\mu(x)=|x-z|^2-|x_0-z|^2$ verify $|\mu|<\eta$ if $|x_0|$ is small enough. We can finish the proof of the claim by choosing $z$ such that the set
	 \[
	 U\cap (B_{|x_0-z|}(z))^c\sim U\cap\{x: \langle x-x_0,e_N\rangle\le 0 \}.
	 \] 
	 This ends the proof.
\end{proof}

\subsection{The limit PDE}
Now, let us show that the half-relaxed limits are sub and supersolutions to our elliptic problem.

\begin{theorem} \label{11}
	The functions $\overline{u}$ and $\underline{u}$  defined in  \eqref{subsuperviscositysols} are respectively viscosity subsolution and supersolutions to the elliptic problem \eqref{PDEPb}.
\end{theorem}

\begin{proof} First, remark that
$\overline{u}$ is upper semicontinuous and $\underline{u}$ is lower semicontinuous.

	Let us write the Dynamic Programming Principle of our parabolic game as
	\[0= \sup_{A} \inf_{B}  \left\{ \frac{1}{\sigma(A\cap B)}
\int_{A\cap B}  \big( u^\varepsilon \big(x + v\varepsilon \big) -u^\varepsilon(x)
\big) 
d \sigma (v)\right\} + K \varepsilon^2. \]
	The use of this equation is the key to pass to the limit
	in the viscosity sense.
	
	Let us show that $\overline{u}$ is a viscosity subsolution to the problem~\eqref{PDEPb}. 
	To this end, 
	let $\phi \in C^{2}$ be a function such as  $\phi-\overline{u}$ achieves a strict minimum at $x_{0}$ in $\overline{B}_{R_0}(x_0)$ for some $R_0>0$ such that $\overline{B}_{R_0}(x_0)$ is in the interior of $\Omega_0$.

Assume first that $\nabla \phi (x_0)\neq 0$. We want to prove that 
$$
\Delta \phi (x_0) \mid_{\langle \nabla \phi (x), v \rangle=0}
 + 1 \geq 0.
$$

	For every $\varepsilon>0$, we can choose $x_\varepsilon \in \overline{B}_{R_0}(x_0)$ such that 
	\[\phi(x_\varepsilon) - u^\varepsilon(x_\varepsilon) - \varepsilon^3\leq \inf_{\overline{B}_{R_0}(x_0)}\{\phi(x)-u^\varepsilon(x)\} \leq \phi(x) - u^\varepsilon(x), \]
	for every $x \in \overline{B}_{R_0}(x_0)$.
	As $\overline{B}_{R_0}(x_0)$ is a compact set, $x_\varepsilon$ converges, up to a subsequence, to some $\overline{x}$ as $\varepsilon \to 0^+$. We have that $\overline{x}= x_0$, since $x_0$ is the only maximum of $\overline{u}-\phi$ in $B_{R_0}(x_0)$. Rearranging the previous expression, we get
	\[u^\varepsilon(y)  - u^\varepsilon(x_\varepsilon) \leq \phi(y) - \phi(x_\varepsilon)  + \varepsilon^3,  \]
	for every $y \in \overline{B}_{R_0}(x_0)$.
	
	Evaluating the DPP at $x = x_\varepsilon$ we obtain 
	\[0\leq \sup_{A} \inf_{B}  \left\{ \frac{1}{\sigma(A\cap B)}
\int_{A\cap B}  \big( \phi \big(x_\varepsilon  + v\varepsilon \big) -\phi (x)
\big) 
d \sigma (v)\right\} + K \varepsilon^2. \]

Since $\phi \in C^2$ performing a second order Taylor expansion and
dividing by $\varepsilon^2$ we obtain 
\begin{equation}
	\label{DPP-44-sect} 0 \!\leq \!   
 \sup_{A} \inf_{B}  \left\{ \frac{1}{\sigma(A\cap B)}
\int_{A\cap B} \Big( \frac{1}{\varepsilon} \langle \nabla \phi (x_\varepsilon), v \rangle +
\frac12
\langle D^2 \phi (x_\varepsilon) v,v \rangle  + o(1) \Big)
d \sigma (v)\right\} + K .
\end{equation}
Now, since we assumed that 
$\nabla \phi (x_0) \neq 0$ and $x_\varepsilon \to x_0$ we also have $\nabla \phi (x_\varepsilon) \neq 0$. Then, after a rotation we can assume that
$\nabla \phi (x_\varepsilon)  = c e_N$ for $\eps$ small enough and $\nabla \phi (x_0)  = c e_N$.
Next we choose $B^*=\{v: \langle e_N, v \rangle \leq \theta_\varepsilon \}$
and we arrive to 
\begin{equation}
	\label{DPP-55-sect} 0 \! \leq \!   
 \sup_{A}   \left\{ \frac{1}{\sigma(A\cap B^*)}
\int_{A\cap B^*} \Big(\frac{1}{\varepsilon} \langle \nabla \phi (x_\varepsilon), v \rangle +
\frac12
\langle D^2 \phi (x_\varepsilon) v,v \rangle  + o(1)\Big)
d \sigma (v)\right\} + K .
\end{equation}
From this inequality we get that for every $A$ with $\sigma (A) \geq 
\frac12 \sigma (S^{N-1}) + \eta_\varepsilon$, we have
$$
 0\geq 
 \frac{1}{\sigma(A\cap B^*)}
\int_{A\cap B^*} \frac{1}{\varepsilon} \langle \nabla \phi (x_\varepsilon), v \rangle d \sigma (v) \geq -C
$$
for some constant $C$ independent of $\varepsilon$.

Using the first inequality, we get
\begin{equation}
	\label{DPP-55-sect-88} 0 \! \leq \!   
 \sup_{A}   \left\{ \frac{1}{\sigma(A\cap B^*)}
\int_{A\cap B^*}
\Big(\frac12
\langle D^2 \phi (x_\varepsilon) v,v \rangle  + o(1)\Big)
d \sigma (v)\right\} + K .
\end{equation}
From the geometric measure lemma (Lemma \ref{lema.clave}) we obtain that
$$
\frac{1}{\sigma(A\cap B^*)}
\int_{A\cap B^*}
\frac12
\langle D^2 \phi (x_\varepsilon) v,v \rangle d \sigma (v)
\to \frac{1}{\mu (\langle e_N,v \rangle=0)}
\int_{\{\langle e_N,v \rangle=0\} }
\frac12
\langle D^2 \phi (x_0) v,v \rangle  d \mu (v).
$$
Hence, passing yo the limit as $\varepsilon \to 0$ we arrive to 
\begin{equation}
	\label{DPP-57-sect} 0 \! \leq \!   
  \left\{  \frac{1}{\mu (\langle e_N,v \rangle=0)}
\int_{\{\langle e_N,v \rangle=0\} }
\frac12
\langle D^2 \phi (x_0) v,v \rangle 
d \mu (v)\right\} + K .
\end{equation}
Now, we just observe that, 
$$
\begin{array}{l}
\displaystyle
 \frac{1}{\mu (\langle e_N, v \rangle=0)}
\int_{\langle e_N, v \rangle=0} \frac12
\langle D^2 \phi (x) v,v \rangle 
d \mu (v) \\[10pt]
\displaystyle
\qquad =\sum_{i,j=1}^{N-1} \frac{1}{\mu (\langle e_N, v \rangle=0)}
\int_{\langle e_N, v \rangle=0} \frac12 \frac{\partial^2\phi}{\partial v_i\partial v_j} (x) v_i,v_j   d\mu (v) \\[10pt]
\displaystyle \qquad = C \sum_{i=1}^{N-1}
\frac{\partial^2 \phi}{\partial v_i^2} (x)
\end{array}
$$
with 
$$
C= \frac12 \vint_{\langle e_N, v \rangle=0} (v_1)^2 \, d \mu (v).
$$

Then, using that we assumed that $\nabla \phi (x_0)$ points in the direction of $e_N$ and that we chose $K=C$ we arrive to
\begin{equation}
	\label{DPP-66-sect} 0 \! \leq \!   
\Delta \phi (x_0) \mid_{\langle \nabla \phi (x_0), v \rangle=0}
 + 1,
\end{equation}
as we wanted to show.

	Thus, since  $\overline{u}(x) \leq 0$ 
for $\partial \Omega_0$ (see Proposition~\ref{ppDatoInicial}), we obtain that $\overline{u}$ is  viscosity subsolution to~\eqref{PDEPb}.

Now, we assume  that $\nabla \phi (x_0)= 0$. We want to prove that 
$$
\Delta \phi (x_0) - \langle  D^2  \phi (x_0) \eta, \eta 
\rangle \geq -1,
$$
for some vector $\eta$ with $|\eta|\leq 1$.

As before, we have that for every there exists a sequence $x_\varepsilon \in \overline{B}_{R_0}(x_0)$ such that 
	\[\phi(x_\varepsilon) - u^\varepsilon(x_\varepsilon) - \varepsilon^3\leq \inf_{\overline{B}_{R_0}(x_0)}\{\phi(x)-u^\varepsilon(x)\} \leq \phi(x) - u^\varepsilon(x), \]
	and  $x_\varepsilon$ converges to $ x_0$. 
Then, as before, we arrive to 
	\[0\leq \sup_{A} \inf_{B}  \left\{ \frac{1}{\sigma(A\cap B)}
\int_{A\cap B}  \big( \phi \big(x_\varepsilon  + v\varepsilon \big) -\phi (x)
\big) 
d \sigma (v)\right\} + k \varepsilon^2.\]

Since $\phi \in C^2$ performing a second order Taylor expansion and
dividing by $\varepsilon^2$ we obtain 
\begin{equation}
	\label{DPP-44-sect-77} 0 \!\leq \!   
 \sup_{A} \inf_{B}  \left\{ \frac{1}{\sigma(A\cap B)}
\int_{A\cap B} \frac{1}{\varepsilon} \langle \nabla \phi (x_\varepsilon), v \rangle +
\frac12
\langle D^2 \phi (x_\varepsilon) v,v \rangle ) + o(1)
d \sigma (v)\right\} + K .
\end{equation}

Now,  we can argue as follows, given $A$ we choose $B$ as 
the symmetrized set of $A$, $B=-A$. Therefore, we have that
$A \cap B = A\cap (-A)$ is symmetric and hence 
$$
 \frac{1}{\sigma(A\cap B)}
\int_{A\cap B} \frac{1}{\varepsilon} \langle \nabla \phi (x_\varepsilon), v \rangle=0
$$
and 
$$
 \frac{1}{\sigma(A\cap B)}
\int_{A\cap B} 
\frac12
\langle D^2 \phi (x_\varepsilon) v,v \rangle ) 
d \sigma (v)= \sum_{i=1}^{N-1} 
\frac{1}{\sigma(A\cap B)}
\int_{A\cap B} 
\frac12
\frac{\partial^2 \phi}{\partial v_i^2} (x_\varepsilon) 
d \sigma (v).
$$

Hence, we arrive to 
\begin{equation}
	\label{DPP-55-sect-000} 0 \! \leq \!   
 \sup_{A} \sum_{i=1}^{N-1} 
\frac{1}{\sigma(A\cap B)}
\int_{A\cap B} 
\frac12
\frac{\partial^2 \phi}{\partial v_i^2} (x_\varepsilon) 
d \sigma (v) + K,
\end{equation}
and then we obtain 
$$
\Delta \phi (x_0) - \langle  D^2  \phi (x_0) \eta, \eta 
\rangle  \rangle  \geq -1.
$$

	To show that $\underline{u}$ is a supersolution, the structure of the proof remains basically the same (with the appropriate changes in the inequalities).

	Thus, since  $\underline{u}(x) \geq 0$ 
for $\partial \Omega_0$ (see Proposition~\ref{ppDatoInicial}), we obtain that $\underline{u}$ is a viscosity supersolution to \eqref{PDEPb}.
\end{proof}

Our aim is to establish the convergence of the value functions of the game, $u^{\varepsilon}$, to the solution to the elliptic problem \eqref{PDEPb}. Having proved that the half relaxed limits satisfy that $\overline{u}$ is a viscosity subsolution and $\underline u$ is a viscosity supersolution to \eqref{PDEPb}, our next goal is to prove that   $\overline u$ and $\underline u$ coincide.  We will achieve this using the Comparison Principle, Theorem \ref{comp.pp}.

\begin{proof}[Proof of Theorem \ref{maintheorem.juego.conv} ] 
	Consider the upper semicontinuous function $\overline{u}$ and  lower semicontinuous function $\underline{u}$ defined by \eqref{subsuperviscositysols}. 
	
	On the one hand, by definition, we have that $\underline{u}(x)\leq \overline{u}(x)$  for all $x\in\Omega_0$.
	
	By Proposition \ref{11}, $\overline{u}$ is subsolution and $\underline{u}$ is a supersolution of problem \eqref{PDEPb}.  Then, $\underline{u}\geq \overline{u}$  in $ \Omega_0$ due to the Comparison Principle, Theorem \ref{comp.pp}. Therefore, we 
	conclude that $$\overline{u}\equiv\underline{u}$$ which implies that $u:=\overline{u}=\underline{u}$ is a continuous function in $ \overline{\Omega}_0$. Moreover, $u$ is the unique solution to the elliptic problem \eqref{PDEPb} and $u^{\varepsilon}$ converges to $u$ pointwise in $\overline{\Omega}_0$. In addition, from these convergences we get that the convergence is uniform $\overline{\Omega}_0$,
	see \cite{Nosotros}. 
	\end{proof}

\subsection{The positivity sets}
Finally, recall that we want to study the positivity sets 
of $u^{\varepsilon}$ and $u$,
\begin{equation}
	\label{OmegaSets.22}	\Omega^{\varepsilon}_{t}:=\{x\,:\;\;u^{\varepsilon}(x)>t\}\;\;\textrm{and}\;\;\Omega_{t}:=\{x\,:\;\;u(x)>t\}.
	\end{equation}
 Next,  as a consequence of the convergence of $u^{\varepsilon}$ to $u$, we will prove that for each fixed $t>0$,
 	\begin{equation*}
 	\Omega_t\subset \liminf_{\varepsilon \to 0} \Omega^\varepsilon_t \subset \limsup_{\varepsilon \to 0} \Omega^\varepsilon_t\subset \overline{\Omega}_t.
 \end{equation*}

\begin{proof}[Proof of Corollary \ref{cor1}]
	Fix $t>0$. We start proving that $$\Omega_t\subset \liminf_{\varepsilon \to 0} \Omega^\varepsilon_t.$$ Suppose that $x\in \Omega_t $. Since 
	\begin{equation*}
		\liminf_{\varepsilon \to 0} \Omega^\varepsilon_t=\bigcup_{\varepsilon_{0}>0}\bigcap_{0<\varepsilon\leq \varepsilon_{0}} \Omega^{\varepsilon}_{t},
	\end{equation*}
	we want to prove that there exists $\varepsilon_{0}>0$ such that $x\in \Omega^{\varepsilon}_{t}$ for all $0<\varepsilon\leq\varepsilon_{0}.$ 
	
	Using that $x\in\Omega_{t}$, we get that there is $\mu>0$ such that $u(x)\geq t+ \mu. $ By Theorem \ref{maintheorem.juego.conv},   $u^{\varepsilon}\to u$. Then, there is $\varepsilon_{0}>0$ such that $u^{\varepsilon}(x) \geq t+\mu/2$ for each $0<\varepsilon\leq \varepsilon_{0}$. Thus, $x\in  \Omega^{\varepsilon}_{t}$ for all $0<\varepsilon\leq\varepsilon_{0}$.

	On the other hand,  suppose that  $x\not\in \overline{\Omega_{t}}$. Our aim is to prove that   $$x\not \in \limsup_{\varepsilon\to 0}\Omega^{\varepsilon}_{t}.$$ Since 
	\begin{equation*}
		\limsup_{\varepsilon \to 0} \Omega^\varepsilon_t=\bigcap_{\varepsilon_{0}>0}\bigcup_{0<\varepsilon\leq \varepsilon_{0}} \Omega^{\varepsilon}_{t},
	\end{equation*} the above is equivalent to prove that there exists $\varepsilon_{0}>0$ for which $x_{0}\not\in \Omega^{\varepsilon}_{t}$ for all $0<\varepsilon\leq \varepsilon_{0}$.
	Due to the fact that $x\not\in\overline{\Omega_{t}}$, we get that there exists $\mu>0$ such that $u(x,t)\leq-\mu$. Using again the convergence of $u^{\varepsilon}$ to $u$, we obtain that there exists $\varepsilon_{0}>0$ such that $u^{\varepsilon}(x)\leq t-\mu/2$ for every $0<\varepsilon\leq \varepsilon_{0}$, that is, for each $0<\varepsilon\leq \varepsilon_{0}$,  $x\not\in\Omega^{\varepsilon}_{t}$.
\end{proof}

\subsection*{Acknowledgments}
	
I. Gonzálvez and J. Ruiz-Cases were supported by grants RED2022-134784-T and CEX2023-001347-S. I. Gonzálvez was also supported by grant PID2019-110712GB-I00. J. Ruiz-Cases was also supported by PID2023-146931NB-I00, all of them funded by MICIU/AEI/10.13039/501100011033 (Spain). A. Miranda and J. D. Rossi were partially supported by CONICET PIP GI No 11220150100036CO (Argentina) and by UBACyT 20020160100155BA (Argentina).

%%% REFERENCES %%%
\EditInfo{May 28, 2025}{June 23, 2025}{Jaqueline Godoy Mesquita, Mariel Sáez Trumper, Rafael Potrie and Tiago Macedo}

\end{document}